\documentclass[12pt, reqno]{amsart}
\usepackage{amsmath, amsthm, amscd, amsfonts, amssymb, graphicx, color}
\usepackage[bookmarksnumbered, colorlinks, plainpages]{hyperref}

\newtheorem{theorem}{Theorem}[section]
\newtheorem{lemma}[theorem]{Lemma}
\newtheorem{proposition}[theorem]{Proposition}
\newtheorem{corollary}[theorem]{Corollary}
\theoremstyle{definition}
\newtheorem{definition}[theorem]{Definition}
\newtheorem{example}[theorem]{Example}

\theoremstyle{remark}
\newtheorem{remark}[theorem]{Remark}
\numberwithin{equation}{section}

\begin{document}
\setcounter{page}{1}

%------------------------------------------------------------------------------------%
%%Don not change any thing in this part

\vspace{1cm}

%------------------------------------------------------------------------------------%

\title[$\Delta$-weak $\phi$-amenability]{on $\Delta$-weak $\phi$-amenability of Banach algebras}

\author[J. laali, M.fozouni]{Javad Laali and Mohammad Fozouni}

\address{Faculty of Mathematics and Computer Science, Department of Mathematics, Kharazmi University,
599 Taleghani, Tehran, Iran.}
\email{\textcolor[rgb]{0.00,0.00,0.84}{Laali@khu.ac.ir}}
\email{\textcolor[rgb]{0.00,0.00,0.84}{fozouni@khu.ac.ir}}

%\dedicatory{This paper3 is dedicated to Professor ABCD}
\subjclass[2010]{46H05, 46H25, 43A30, 22D15.}
\keywords{Banach algebra, Character, Approximate Identity, Locally Compact Group}
%\date{Received: xxxxxx; Revised: yyyyyy; Accepted: zzzzzz.}
%\newline \indent $^{*}$ Corresponding author}

\begin{abstract}
Let $A$ be a Banach algebra and $\phi\in \Delta(A)\cup\{0\}$. We introduce and study the notion  $\Delta$-weak $\phi$-amenability of Banach algebra $A$. It is shown that $A$ is  $\Delta$-weak $\phi$-amenable if and only if $\ker(\phi)$ has a bounded $\Delta$-weak approximate identity.  We prove that every $\Delta$-weak $\phi$-amenable Banach algebra has a bounded $\Delta$-weak approximate identity. Finally, we  examine this notion for some algebras over locally compact groups and give a characterization of $\Delta$-weak $\phi$-amenablity of Figa-Talamanca-Herz algebras.
\end{abstract}
\maketitle
\section{Introduction}
Let  $A$ be a Banach algebra, $\Delta(A)$ be the character space of $A$, i.e., the space of all nonzero homomorphisms from $A$ into $\mathbb{C}$ and $A^{*}$ be the dual space of $A$ consisting of all bounded linear functional from $A$ into $\mathbb{C}$. %For each character $\phi\in \Delta(A)$, by \cite[Lemma 2.1.5]{Kaniuth2} we know that $||\phi||\leq 1$. Therefore, $\Delta(A)\subseteq A^{*}$.

Throughout this paper we assume that $A$ is a Banach algebra such that $\Delta(A)\neq \emptyset$.

Let $\{e_{\alpha}\}$ be a net in a Banach algebra $A$. The net $\{e_{\alpha}\}$ is called,
\begin{enumerate}
  \item an \emph{approximate identity} if, for each $a\in A$, $||ae_{\alpha}-a||+||e_{\alpha}a-a||\rightarrow 0$,
  \item a \emph{weak approximate identity} if, for each $a\in A$, $|f(ae_{\alpha})-f(a)|+|f(e_{\alpha}a)-f(a)|\rightarrow 0$ for all $f\in A^{*}$,
  \item a \emph{$\Delta$-weak approximate identity} if, for each $a\in A$, $|\phi(e_{\alpha}a)-\phi(a)|\rightarrow 0$ for all $\phi\in \Delta(A)\cdot$
\end{enumerate}
The notion of weak approximate identity was originally introduced for the study of the second dual $A^{**}$ of a Banach algebra $A$. For technical reasons, bounded approximate identities are of interest for mathematicians. It is proved that every Banach algebra $A$ which has a bounded weak approximate identity, also  has a bounded approximate identity and  conversely (\cite[Proposition 33.2]{Doran}).  But  in \cite{Jones}, Jones and Lahr proved that the approximate identity and $\Delta$-weak approximate identity of a Banach algebra are different. They showed that there exists a Banach algebra $A$ for which has a bounded $\Delta$-weak approximate identity, but it does  not have any approximate identity. Indeed, if $S=\mathbb{Q}^{+}$ is the semigroup of positive rationales under addition, they showed that the semigroup algebra $l^{1}(S)$ has a bounded $\Delta$-weak approximate identity, but it does not have any bounded or unbounded approximate identity.

%For the convenience of the reader, here we give some of the important definitions and results from \cite{Jones} and \cite{Kaniuth}. We generally follow \cite{Dales} for some of notations and concepts.
\begin{definition}\label{maindef}
Let $A$ be a Banach algebra. A $\Delta$-weak approximate identity for subspace $B\subseteq A$ is a net $\{a_{\alpha}\}$ in $B$ such that
\begin{equation*}
\lim_{\alpha}|\phi(aa_{\alpha})-\phi(a)|=0\hspace{0.5cm}(a\in B, \phi\in \Delta(A))\cdot
\end{equation*}
\end{definition}

For $a\in A$ and $f\in A^{*}$,  the linear functional $f.a$ defined as follows
 \begin{center}
 $<f.a,b>=<f,ab>=f(ab)\hspace{0.5cm}  (b\in A)\cdot$
 \end{center}
\begin{definition}
Let $A$ be a Banach algebra and $\phi\in \Delta(A)$. The Banach algebra $A$ is said to be \emph{$\phi$-amenable} if, there exists an $m\in A^{**}$ such that satisfies the following relations,
\begin{enumerate}
  \item $m(\phi)=1$,
  \item $m(f.a)=\phi(a)m(f)\hspace{0.5cm}(a\in A, f\in A^{*})\cdot$
\end{enumerate}
\end{definition}

The concept of character amenability first introduced by Monfared in \cite{Monfared}. Also, Kaniuth, Lau and Pym in \cite{Kaniuth}, investigated the concept of $\phi$-amenability of Banach algebras and give the following result.

\begin{theorem}\label{main} Let A be a Banach algebra and $\phi\in \Delta(A)$. Then $\ker(\phi)$ has a bounded right approximate identity if and only if $A$ is $\phi$-amenable and has a bounded right approximate identity.
\end{theorem}
\begin{proof}
See \cite[Corollary 2.3]{Kaniuth}.
\end{proof}
In the next section of this paper, first we give the basic definition of our work, that is $\Delta$-weak $\phi$-amenability of a Banach algebra $A$. Then characterize it through the existence of a bounded $\Delta$-weak approximate identity for $\ker(\phi)$. Suppose that $A$ is a Banach algebra and $I$ is a closed ideal of $A$ such that $I$ and $A/I$ both have  bounded approximate identities. Then $A$ has a bounded approximate identity (\cite[Proposition 7.1]{Doran}). We give a variant of  this theorem. With use of this theorem, we prove that each Banach algebra $A$ that is $\Delta$-weak $\phi$-amenable,  has a bounded $\Delta$-weak approximate identity.

In Section 3, we investigate some of the hereditary  properties of $\Delta$-weak $\phi$-amenablity. In section 4 we study group algebras and Figa-Talamanca Herz algebras of a locally compact group with respect to this notion and prove that when $1<p<\infty$ and $\phi\in \Delta(A_{p}(G)\cup\{0\}$, $A_{p}(G)$ is $\Delta$-weak $\phi$-amenable if and only if $G$ is an amenable group.

In the final section, we just give some examples which shows the different situations that might be happen for  definitions  and theorems in Section 2 and 3.

%%%%%%%%%%%%%%%%%%%%%%%%%%%
% Next section
%%%%%%%%%%%%%%%%%%%%%%%%%%%

\section{Main definition and its characterization}
We commence this section with the following definition that is our main concern.
\begin{definition}
Let $A$ be a Banach algebra and $\phi\in \Delta(A)\cup \{0\}$. We say that $A$ is \emph{$\Delta$-weak $\phi$-amenable} if, there exists an $m\in A^{**}$ such that $m(\phi)=0$ and $m(\psi.a)=\psi(a)$ for each $a\in \ker(\phi)$ and $\psi\in \Delta(A)$.
\end{definition}

%%%%%%%%%%%%%%%%%%%%%%%%%%%%%%%%%%%%%%%%%%%%%%%%%%%%%%%%%%%%%%%%%%%%%%%
%%%%%%%%%%%%%%%%%%%%%%%%%%%%%%%%%%%%%%%%%%%%%%%%%%%%%%%%%%%%%%%%%%%%%%%
%main theorem
%%%%%%%%%%%%%%%%%%%%%%%%%%%%%%%%%%%%%%%%%%%%%%%%%%%%%%%%%%%%%%%%%%%%%%%
%%%%%%%%%%%%%%%%%%%%%%%%%%%%%%%%%%%%%%%%%%%%%%%%%%%%%%%%%%%%%%%%%%%%%%%

Now, we characterize the concept of $\Delta$-weak $\phi$-amenability as follows. Recall that if $A$ is a Banach algebra, for each $a\in A$, $\widehat{a}\in A^{**}$ defined by $\widehat{a}(f)=f(a)$ for all $f\in A^{*}$.
%===================================================================
\begin{theorem}\label{T: main charac} Let $A$ be a Banach algebra and $\phi\in \Delta(A)\cup \{0\}$. Then $A$ is $\Delta$-weak $\phi$-amenable if and only if $\ker(\phi)$ has a bounded $\Delta$-weak approximate identity.
\end{theorem}
\begin{proof}
Let $\{e_{\alpha}\}$ be a bounded $\Delta$-weak approximate identity for $\ker(\phi)$. So, $\{\widehat{e_{\alpha}}\}$ is a bounded net in $A^{**}$. Therefore, the Banach-Alaoglu's Theorem  (\cite[Theorem A.3.20]{Dales}) yields $\{\widehat{e_{\alpha}}\}$ is a compact set with $w^{*}$-topology and hence has a $w^{*}$-accumulation point $m$, i.e., there exists a subnet that we show it  again by $\{\widehat{e_{\alpha}}\}$ such that $m=w^{*}-\lim_{\alpha}(\widehat{e_{\alpha}})$.

Now,  we have
\begin{center}
$m(\phi)=\lim_{\alpha}\widehat{e_{\alpha}}(\phi)=\lim_{\alpha}\phi(e_{\alpha})=0,$
\end{center}
and for all $\psi\in \Delta(A)$ and $a\in \ker(\phi)$,
\begin{center}    $m(\psi.a)=\lim_{\alpha}\widehat{e_{\alpha}}(\psi.a)=\lim_{\alpha}\psi.a(e_{\alpha})=\lim_{\alpha}\psi(ae_{\alpha})=\psi(a)\cdot$
\end{center}

Conversely, since $m\in A^{**}$, Goldstine's Theorem (\cite[Theorem A.3.29(i)]{Dales}) yields  there exists a net $\{e_{\alpha}\}$ in $A$ such that $m=w^{*}-\lim_{\alpha}\widehat{e_{\alpha}}$ and $||e_{\alpha}||\leq ||m||$. So, $\{e_{\alpha}\}$ is a bounded net such that for all $\psi\in \Delta(A)$ and $a\in \ker(\phi)$ we have
\begin{align*}
\psi(a)=m(\psi.a)&=\lim_{\alpha}\widehat{e_{\alpha}}(\psi.a)\\
&=\lim_{\alpha}\psi.a(e_{\alpha})\\
&=\lim_{\alpha}\psi(ae_{\alpha})\cdot
\end{align*}
Hence $\lim_{\alpha}|\psi(ae_{\alpha})-\psi(a)|=0$.

Note that the net $\{e_{\alpha}\}$ is not a subset of $\ker(\phi)$. To construct a net in $\ker(\phi)$,  we have two cases $\phi=0$ or $\phi\in \Delta(A)$. If $\phi=0$, then the net $\{e_{\alpha}\}$ is a bounded $\Delta$-weak approximate identity for $A=\ker(\phi)$.

If  $\phi\in \Delta(A)$, then there exists $x_{0}\in A$ such that $\phi(x_{0})\neq 0$. Put $a_{0}=\frac{x_{0}}{\phi(x_{0})}$ and suppose that $a_{\alpha}=e_{\alpha}-\phi(e_{\alpha})a_{0}$ for all $\alpha$. Obviously $\{a_{\alpha}\}\subseteq \ker(\phi)$. On the other hand, since $m(\phi)=0$, we conclude that $\lim_{\alpha}\phi(e_{\alpha})=0$. Hence
\begin{equation*}
\lim_{\alpha}|\psi(aa_{\alpha})-\psi(a)|=0\hspace{0.5cm}(a\in \ker(\phi), \psi\in \Delta(A))\cdot
\end{equation*}
\end{proof}
%%%%%%%%%%%%%%%%%%%%%%%%%%%%%%%%%%%%%%%%%%%%%%%%%%%%%%%%%%%%%%%%%%%%%%%
For simplicity of notation, we let b.$\Delta$-w.a.i stand for bounded approximate identity and b.a.i stand for bounded approximate identity.

By Theorems \ref{main} and \ref{T: main charac} one can see that every $\phi$-amenable Banach algebra which has a bounded right approximate identity is $\Delta$-weak $\phi$-amenable. But the converse of this assertion is not valid in general.
%%%%%%%%%%%%%%%%%%%%%%%%%%%%%%%%%%%%%%%%%%%%%%%%%%%%%%%%%%%%%%%%%%%%%%%
%%%%%%%%%%%%%%%%%%%%%%%%%%%%%%%%%%%%%%%%%%%%%%%%%%%%%%%%%%%%%%%%%%%%%%%
%%%%%%%%Remark%%%%%%%%%%%%%%%%%%%%%%%%%%%%%%%%%%%%%%%%%%%%%%%%%%%%%%%%
\begin{remark} Let $A$ be a $\Delta$-weak $\phi$-amenable Banach algebra and $\{e_{\alpha}\}$ be a $\Delta$-weak approximate identity of $\ker(\phi)$. If there exists $a_{0}\in A$ with $\phi(a_{0})=1$ and $\lim_{\alpha}|\psi(a_{0}e_{\alpha})-\psi(a_{0})|=0$ for all $\psi\in \Delta(A)\setminus \{\phi\}$, then there exists $m\in A^{**}$ such that
\begin{enumerate}\label{E: 1}
  \item $m(\phi)=0,$
  \item $m(\psi.a)=\psi(a)\hspace{0.5cm}(a\in A, \psi\in \Delta(A)\setminus\{\phi\})\cdot$
\end{enumerate}
By a similar argument as in the above theorem we can show the existence of $m$.  Let $a\in A$. It is clear that $a-\phi(a)a_{0}\in \ker(\phi)$. So, for all $\psi\in \Delta(A)\setminus\{\phi\}$ we have $m(\psi.(a-\phi(a)a_{0}))=\psi(a-\phi(a)a_{0})$. Therefore $m(\psi.a)=\psi(a)$.

Note that in part (2), it is necessary that $\psi\neq \phi$. If $\psi=\phi$ we have $m(\phi.a)=\phi(a)$. On the other hand, there exists an $a\in A$ such that $\phi(a)\neq 0$. So, we have
\begin{equation*}
\phi(a)=m(\phi.a)=m(\phi(a)\phi)=\phi(a)m(\phi)\cdot
\end{equation*}
Therefore, $m(\phi)=1$ and this is a contradiction.
\end{remark}
%%%%%%%%%%%%%%%%%%%%%%%%%%%%%%%%%%%%%%%%%%%%%%%%%%%%%%%%%%%%%%%%%%%%%%%
%%%%%%%%%%%%%%%%%%%%%%%%%%%%%%%%%%%%%%%%%%%%%%%%%%%%%%%%%%%%%%%%%%%%%%%
%%%%%%%%%%%%%%%%%%%%%%%%%%%%%%%%%%%%%%%%%%%%%%%%%%%%%%%%%%%%%%%%%%%%%%%
The following lemma is needed in the sequel.
\begin{lemma}\label{L: two character}Let $A$ be a Banach algebra such that $\phi,\psi\in \Delta(A)$ and $\phi\neq \psi$. Then there exists $a\in A$ such that $\phi(a)=0$ and $\psi(a)=1$.
 \end{lemma}
\begin{proof}
See the proof of \cite[Theorem 3.3.14]{Kaniuth2}.
\end{proof}
%===========================================================
\begin{remark}\label{remark}
If $A$ is a Banach algebra with $\Delta(A)\setminus\{\phi\}\neq\emptyset$, then $A$ is $\Delta$-weak $\phi$-amenable if and only if there exists a bounded net $\{e_{\alpha}\}$ in $\ker(\phi)$ such that  $\lim_{\alpha}\psi(e_{\alpha})=1$ for each $\psi\in \Delta(A)\setminus\{\phi\}$ or equivalently there exists an $m\in A^{**}$ with $m(\phi)=0$ and $m(\psi)=1$ for all $\psi\in \Delta(A)\setminus\{\phi\}$.
\end{remark}

%===============================================================================
The following proposition allow us to produce Banach algebras which are $\Delta$-weak $\phi$-amenable, but they are not $\phi$-amenable.
\begin{proposition}\label{mainprop} Let $A$ be a Banach algebra such that $0<|\Delta(A)|\leq2$, i.e., $\Delta(A)$ has 1 or 2 elements. Then $A$ is $\Delta$-weak $\phi$-amenable.
\end{proposition}
\begin{proof} If $A$ only has one character, the proof is easy. Therefore, we omit it. In the second case, let $\Delta(A)=\{\phi, \psi\}$ and $\phi\neq \psi$. Hence, by Lemma \ref{L: two character} there exists an $e\in A$ with $\phi(e)=0$ and $\psi(e)=1$. Now, put $m=\widehat{e}$. Clearly, $m(\phi)=0$ and $m(\psi)=1$, so by  Remark \ref{remark} $A$ is $\Delta$-weak $\phi$-amenable.
\end{proof}
%%%%%%%%%%%%%%%%%%%%%%%%%%%%%%%%%%%%%%%%%%%%%%%%%%%%%%%%%%%%%%%%%%%%%%%
%%%%%%%%%%%%%%%%%%%%%%%%%%%%%%%%%%%%%%%%%%%%%%%%%%%%%%%%%%%%%%%%%%%%%%%
%%%%%%%%%%%%%%%%%%%%%%%%%%%%%%%%%%%%%%%%%%%%%%%%%%%%%%%%%%%%%%%%%%%%%%%
The following theorem is a useful tools in the rest of this section.
\begin{theorem}\label{T:Ap identity} Let $A$ be a Banach algebra and $I$ be a closed two-sided  ideal of $A$ which has a b.$\Delta$-w.a.i and the quotient Banach algebra $A/I$ has a  bounded left approximate identity (b.l.a.i). Then $A$ has a b.$\Delta$-w.a.i.
\end{theorem}
\begin{proof} Let $\{e_{\alpha}\}$ be a b.$\Delta$-w.a.i for $I$ and $\{f_{\delta}+I\}$ be a b.l.a.i for $A/I$. Suppose that $F=\{a_{1},...a_{m}\}$ is a finite subset of $A$ and $n$ is a positive integer. Let $M$ be an upper bound for $\{||e_{\alpha}||\}$. For $\lambda=(F,n)$, there exists $f_{\delta_{\lambda}}$ such that
\begin{equation*}
||f_{\delta_{\lambda}}a_{i}-a_{i}+I||<\frac{1}{2(1+M)n}\hspace{0.5cm}(i=1,2,3,...,m)\cdot
\end{equation*}
Therefore, there exists $y_{i}\in I$ such that
\begin{equation*}
||f_{\delta_{\lambda}}a_{i}-a_{i}+y_{i}||<\frac{1}{2(1+M)n}\hspace{0.5cm}(i=1,2,3,...,m)\cdot
\end{equation*}
Let $\psi\in \Delta(A)$. Since $\{e_{\alpha}\}$ is a b.$\Delta$-w.a.i, for each $y_{i}$ with $i\in\{1,2,3,...,m\}$ which satisfy the above relation, there exists $e_{\alpha_{\lambda}}\in\{e_{\alpha}\}$ such that
\begin{equation*}
||\psi(e_{\alpha_{\lambda}}y_{i})-\psi(y_{i})||<\frac{1}{2n}\hspace{0.5cm}(i=1,2,3,...,m)\cdot
\end{equation*}
Now, for each $i\in\{1,2,3,...,m\}$ we have
\begin{align*}
||\psi((e_{\alpha_{\lambda}}+f_{\delta_{\lambda}}-e_{\alpha_{\lambda}}f_{\delta_{\lambda}})a_{i})-\psi(a_{i})||&\leq ||\psi(f_{\delta_{\lambda}}a_{i}-a_{i}+y_{i})||\\
&+||\psi(e_{\alpha_{\lambda}}y_{i})-\psi(y_{i})||\\
&+||\psi(e_{\alpha_{\lambda}}a_{i}-e_{\alpha_{\lambda}}f_{\delta_{\lambda}}a_{i}-e_{\alpha_{\lambda}}y_{i})||\\
&\leq ||f_{\delta_{\lambda}}a_{i}-a_{i}+y_{i}||\\
&+\frac{1}{2n}+M||a_{i}-f_{\delta_{\lambda}}a_{i}-y_{i}||\\
&<\frac{1}{n}\cdot
\end{align*}
Therefore, $\{e_{\alpha_{\lambda}}+f_{\delta_{\lambda}}-e_{\alpha_{\lambda}}f_{\delta_{\lambda}}\}_{\lambda\in \Lambda}$ is a $\Delta$-w.a.i for $A$, where $\Lambda=\{(F,n):F\subseteq A \textmd{ is  finite}, n\in \mathbb{N}\}$ is a directed set with $(F_{1},n_{1})\leq (F_{2},n_{2})$ if $F_{1}\subseteq F_{2}$ and $n_{1}\leq n_{2}$.

Now, we show that there exists a b.$\Delta$-w.a.i for $A$. Since $\{f_{\delta}+I\}$ is bounded, there exists a positive integer $K$ such that $||f_{\delta}+I||<K$ for each $\delta$. So,  there exists $y_{\delta}\in I$ such that $||f_{\delta}+I||<||f_{\delta}+y_{\delta}||<K$. Put $f^{'}_{\delta}=f_{\delta}+y_{\delta}$. Hence, $\{f^{'}_{\delta}+I\}$ is a bounded approximate identity for $A/I$ which $\{f^{'}_{\delta}\}$ is  bounded.
Now, we have
\begin{center}
$||e_{\alpha_{\lambda}}+f^{'}_{\delta_{\lambda}}-e_{\alpha_{\lambda}}f^{'}_{\delta_{\lambda}}||\leq ||e_{\alpha_{\lambda}}||+||f^{'}_{\delta_{\lambda}}||+||e_{\alpha_{\lambda}}||||f^{'}_{\delta_{\lambda}}||<M+K+KM\cdot$
\end{center}
Therefore, $A$ has a b.$\Delta$-w.a.i.
\end{proof}
It is straightforward to see that for every closed two-sided ideal $I$ with codimension one of a Banach algebra $A$, the quotient Banach algebra $A/I$ has a bounded approximate identity. So, we have the following corollary.
\begin{corollary}\label{corollary}
Let $A$ be a Banach algebra and $\phi\in \Delta(A)\cup\{0\}$. If $A$ is $\Delta$-weak $\phi$-amenable, then $A$ has a b.$\Delta$-w.a.i.
\end{corollary}
\begin{proof}
Since $A$ is $\Delta$-weak $\phi$-amenable, by Theorem \ref{T: main charac}, $\ker(\phi)$ has a b.$\Delta$-w.a.i. Also, $A/\ker(\phi)$ has a b.a.i, because the codimension of $\ker(\phi)$ is one. Then, by Theorem \ref{T:Ap identity}, $A$ has a b.$\Delta$-w.a.i.
\end{proof}
%============================================================
\begin{corollary}\label{Cor: 0-amen}
Let $A$ be a Banach algebra and $\phi\in \Delta(A)$. If $A$ is $\Delta$-weak $\phi$-amenable, then $A$ is $\Delta$-weak 0-amenable.
\end{corollary}

The converse of the above corollary is not valid in general (see Example \ref{ex2}).

Note that there exists Banach algebras which do not have any b.$\Delta$-w.a.i. As an example, consider $S_{p}(G)=L^{1}(G)\cap L^{p}(G)$ with norm defined by $||f||_{S_{p}(G)}=\max\{||f||_{1},||f||_{p}\}$ which is a Segal algebra (see \cite{Rietr} for a full discussion on Segal algebras). Now, by \cite[Remark 2]{Inoue}, if $G$ is an infinite abelian  compact group, $S_{p}(G)$ has no b.$\Delta$-w.a.i.
%%%%%%%%%%%%%%%%%%%%%%%%%%%%%%%%%%%%%%%%%%%%%%%%%%%%%%%%%%%%%%%%%%%%%%%
%%%%%%%%%%%%%%%%%%%%%%%%%%%%%%%%%%%%%%%%%%%%%%%%%%%%%%%%%%%%%%%%%%%%%%%
%\Delta-weak identity
%%%%%%%%%%%%%%%%%%%%%%%%%%%%%%%%%%%%%%%%%%%%%%%%%%%%%%%%%%%%%%%%%%%%%%%
%%%%%%%%%%%%%%%%%%%%%%%%%%%%%%%%%%%%%%%%%%%%%%%%%%%%%%%%%%%%%%%%%%%%%%%

There exists a $\Delta$-weak version of an identity in a Banach algebra $A$.
\begin{definition} Let $A$ be a Banach algebra. We say that $e\in A$ is a $\Delta$-weak identity for $A$ if for each $\phi\in \Delta(A)$, $\phi(e)=1$ or equivalently
\begin{equation*}
\phi(ea)=\phi(a)\hspace{0.5cm}(a\in A, \phi\in \Delta(A))\cdot
\end{equation*}
\end{definition}
It is obvious that the identity of a Banach algebra $A$ is a $\Delta$-weak identity of $A$, but the converse is not valid in general.
%======================================================================

The following theorem gives a necessary condition for $\Delta$-weak $\phi$-amenability of finite dimensional Banach algebras.
%====================================================================
\begin{theorem} Let $A$ be a finite dimensional Banach algebra. If $A$  is $\Delta$-weak $\phi$-amenable, then it has a $\Delta$-weak identity.
\end{theorem}
\begin{proof}
In view of Corollary \ref{corollary}, $A$ has a b.$\Delta$-w.a.i, say $\{e_{\alpha}\}$. By the Heine-Borel's Theorem (\cite[Theorem 2.38]{Nair}), we know that every closed and bounded subset of a finite dimensional normed linear space is compact.

So, $\{e_{\alpha}\}$ is compact. Therefore, there exists  $e\in A$ and a convergent subnet that we show it again  by $\{e_{\alpha}\}$ such that converges to $e$. Now, for each $\psi\in \Delta(A)$ and $a\in A$, we have
\begin{align*}
\psi(ea)=\lim_{\alpha}\psi(e_{\alpha}a)=\psi(a)\cdot
\end{align*}
Therefore, $e$ is a $\Delta$-weak identity for $A$.
\end{proof}
%%%%%%%%%%%%%%%%%%%%%%%%%%%%%%%%%%%%%%%%%%%%%%%%%%%%%%%%%%%%%%%%%%%%%%%
%%%%%%%%%%%%%%%%%%%%%%%%%%%%%%%%%%%%%%%%%%%%%%%%%%%%%%%%%%%%%%%%%%%%%%%
%Hereditary properties
%%%%%%%%%%%%%%%%%%%%%%%%%%%%%%%%%%%%%%%%%%%%%%%%%%%%%%%%%%%%%%%%%%%%%%%
%%%%%%%%%%%%%%%%%%%%%%%%%%%%%%%%%%%%%%%%%%%%%%%%%%%%%%%%%%%%%%%%%%%%%%%
\section{Hereditary properties}
In this section we give some of the hereditary properties of $\Delta$-weak $\phi$-amenability.

\begin{theorem}\label{T: HP dense range} Let $A$ and $B$ be Banach algebras, $\phi\in \Delta(B)$ and $h:A\rightarrow B$ be a dense range continuous homomorphism. If $A$ is $\Delta$-weak $\phi \circ h$-amenable, then $B$ is $\Delta$-weak $\phi$-amenable.
\end{theorem}
\begin{proof} Let $A$ be $\Delta$-weak $\phi \circ h$-amenable. So, there exists $m\in A^{**}$ such that, $m(\phi \circ h)=0$ and $m(\psi.a)=\psi(a)$ for all $a\in \ker(\phi \circ h)$ and $\psi\in \Delta(A)$.

Define $n\in B^{**}$ as follows
\begin{equation*}
n(g)=m(g\circ h)\hspace{0.5cm}(g\in B^{*})\cdot
\end{equation*}
So, $n(\phi)=m(\phi \circ h)=0$. For each $b\in\ker(\phi)$, there exists a sequence $\{e_{n}\}$ in $A$ such that $\lim_{n}h(e_{n})=b$. Put $a_{n}=e_{n}-\phi \circ h(e_{n})a_{0}$ where $\phi \circ h(a_{0})=1$. It is obvious that $a_{n}\in \ker (\phi \circ h)$ for each $n$ and $\lim_{n}h(a_{n})=b$. Also, for each $\psi^{'}\in \Delta(B)$, $(\psi^{'}\circ h).a_{n}\rightarrow (\psi^{'}.b)\circ h$ in $A^{*}$, since
\begin{align*}
||(\psi^{'}\circ h).a_{n}-(\psi^{'}.b)\circ h||&=\sup_{||a||\leq 1}||\psi^{'}(h(a_{n}a))-\psi^{'}(bh(a))||\\
&\leq \sup_{||a||\leq 1}||h(a_{n})h(a)-bh(a)||\\
&\leq ||h(a_{n})-b||||h||\cdot
\end{align*}
Therefore, for all $\psi^{'}\in \Delta(B)$ we have
\begin{align*}
n(\psi^{'}.b)=m((\psi^{'}.b)\circ h)&=\lim_{n}m((\psi^{'}\circ h).a_{n})\\
&=\lim_{n}\psi^{'} \circ h(a_{n})\\
&=\lim_{n}\psi^{'}(h(a_{n}))\\
&=\psi^{'}(b)\cdot
\end{align*}
So, $B$ is $\Delta$-weak $\phi$-amenable.
\end{proof}
Let $I$ be a closed ideal of a Banach algebra $A$. As mentioned in the proof of \cite[Theorem 2.6 (ii)]{Monfared2} one can see if $I$ has an approximate identity, then every $\phi\in \Delta(I)$ extend to some $\widetilde{\phi}\in \Delta(A)$. To see this let $\{e_{\alpha}\}$ be an approximate identity of $I$ and  $u\in I$ be an element with  $\phi(u)=1$. If $a\in A$ and $b\in \ker(\phi)$, then we have
\begin{equation*}
\phi(ab)=\lim_{\alpha}\phi(ae_{\alpha}b)=0\cdot
\end{equation*}
Therefore, $ab\in \ker(\phi)$ and this shows that $\ker(\phi)$ is a left ideal in $A$. Now, Define $\widetilde{\phi}:A\rightarrow \mathbb{C}$ by $\widetilde{\phi}(a)=\phi(au)$ for all $a\in A$. Since for $a,b\in A$, $bu-ubu\in \ker(\phi)$ therefore, $abu-aubu=a(bu-ubu)\in \ker(\phi)$. So,  we conclude that $\phi(abu)=\phi(au)\phi(bu)$. Hence
\begin{equation*}
\widetilde{\phi}(ab)=\phi(abu)=\phi(au)\phi(bu)=\widetilde{\phi}(a)\widetilde{\phi}(b)\hspace{0.5cm}(a, b\in A)\cdot
\end{equation*}

\begin{proposition}\label{prop1} Let $A$ be a Banach algebra, $I$ be a closed ideal of $A$ such that has a bounded approximate identity and  $\phi\in \Delta(A)$ with $I\nsubseteq \ker(\phi)$. If $A$ is $\Delta$-weak $\phi$-amenable, then $I$ is  $\Delta$-weak $\phi_{|I}$-amenable.
%%%%%%%%%%%%%%%%%%%%%%%%%%%%%%%%%%%%%%%%%%%%%%%%%%%%%%%%%%%%%%%%%%%%%%%
%%%%%%%%%%%%%%%%%%%%%%%%%%%%%%%%%%%%%%%%%%%%%%%%%%%%%%%%%%%%%%%%%%%%%%%
%(Why the ideal I is closed?)
%we say the ideal I is closed because we want to speak of \Delta(I), usually we speak of character space of a Banach space.
%%%%%%%%%%%%%%%%%%%%%%%%%%%%%%%%%%%%%%%%%%%%%%%%%%%%%%%%%%%%%%%%%%%%%%%
%%%%%%%%%%%%%%%%%%%%%%%%%%%%%%%%%%%%%%%%%%%%%%%%%%%%%%%%%%%%%%%%%%%%%%%
\end{proposition}

\begin{proof} Let $\{a_{\beta}\}$ be a bounded approximate identity for $I$ and $\{e_{\alpha}\}$ be a b.$\Delta$-w.a.i for $\ker(\phi)$. It is clear that $\lim_{\beta}\psi(a_{\beta})=1$ for all $\psi\in \Delta(I)$.

Put $c_{(\alpha,\beta)}=e_{\alpha}a_{\beta}$ for all $\alpha, \beta$. Then $\{c_{(\alpha,\beta)}\}_{(\alpha,\beta)}$ is a bounded net in $I$. Now, for each $a\in \ker(\phi_{|I})$ and $\psi\in \Delta(I)$ we have
\begin{align*}
\lim_{(\alpha,\beta)}\psi(ac_{(\alpha,\beta)})&=\lim_{(\alpha,\beta)}\psi(ae_{\alpha}a_{\beta})\\
&=(\lim_{(\alpha,\beta)}\psi(ae_{\alpha}))(\lim_{(\alpha,\beta)}\psi(a_{\beta}))\\
&=(\lim_{(\alpha,\beta)}\widetilde{\psi}(ae_{\alpha}))(\lim_{(\alpha,\beta)}\psi(a_{\beta}))\\
&=\widetilde{\psi}(a)=\psi(a)\cdot
\end{align*}
Therefore, $I$ is $\Delta$-weak $\phi_{|I}$-amenable by Theorem \ref{T: main charac}.
\end{proof}

For each Banach algebra $A$, we can extend each $\phi\in \Delta(A)$ uniquely %why unique?
to a  character $\phi$ of $A^{**}$ defined by $\widehat{\phi}(a^{**})=a^{**}(\phi)$ for all $a^{**}\in A^{**}$. So, we have the following result.

\begin{proposition} Let $A$ be a Banach algebra and $\phi\in \Delta(A)$. If  $A^{**}$ is $\Delta$-weak $\widehat{\phi}$-amenable, then $A$ is $\Delta$-weak $\phi$-amenable.
\end{proposition}
\begin{proof} Let $A^{**}$ be $\Delta$-weak $\widehat{\phi}$-amenable. So, there exists $m^{**}\in A^{****}$ which satisfies the following relations,
\begin{enumerate}
  \item $m^{**}(\widehat{\phi})=0,$
  \item $m^{**}(\Psi.F)=\Psi(F)\hspace{0.5cm}(\Psi\in \Delta(A^{**}), F\in \ker(\widehat{\phi}))\cdot$
\end{enumerate}
Put $m(f)=m^{**}(\widehat{f})$ for all $f\in A^{*}$. Therefore, $m(\phi)=m^{**}(\widehat{\phi})=0$ and for each $a\in \ker(\phi)\subseteq \ker(\widehat{\phi})$ we have
\begin{center}
$m(\psi.a)=m^{**}(\widehat{\psi.a})=m^{**}(\widehat{\psi}.\widehat{a})=\widehat{\psi}(\widehat{a})=\widehat{a}(\psi)=\psi(a)\quad(\psi\in \Delta(A))\cdot$
\end{center}
Therefore, $A$ is $\Delta$-weak $\phi$-amenable.
\end{proof}

%%%%%%%%%%%%%%%%%%%%%%%%%%%%%%%%%
%%%%%%%%%%%%%%%%%%%%%%%%%%%%%%%%%%
% an example of a Banach algebra that is delta-weak \phi-amenable but is not delta weak \widetilde{\phi}-amenable?????
%%%%%%%%%%%%%%%%%%%%%%%%%%%%%%%%%
%%%%%%%%%%%%%%%%%%%%%%%%%%%%%%%%%%
\section{Some results on algebras over locally compact groups}
Let $G$ be a locally compact group. For $1<p<\infty$ let $A_{p}(G)$ denote the subspace of $C_{0}(G)$ consisting of functions of the form $u=\sum_{i=1}^{\infty}f_{i}\ast \widetilde{g_{i}}$ where $f_{i}\in L^{p}(G)$, $g_{i}\in L^{q}(G)$, $1/p+1/q=1$, $\sum_{i=1}^{\infty}||f_{i}||_{p}||g_{i}||_{q}<\infty$ and $\widetilde{f}(x)=\overline{f(x^{-1})}$ for all $x\in G$. $A_{p}(G)$ is called the Figa-Talamanca-Herz algebra and with the pointwise operations and the following norm is a Banach algebra,
\begin{equation*}
||u||_{A_{p}(G)}=\inf\{\sum_{i=1}^{\infty}||f_{i}||_{p}||g_{i}||_{q}: u=\sum_{i=1}^{\infty}f_{i}\ast \widetilde{g_{i}}\}\cdot
\end{equation*}
It is obvious that for each $u\in A_{p}(G)$, $||u||\leq ||u||_{A_{p}(G)}$ where $||u||$ is the norm of $u$ in $C_{0}(G)$. Also we know that $\Delta(A_{p}(G))=G$, i.e., each character of $A_{p}(G)$ is an evaluation function at some $x\in G$ \cite[Theorem 3]{Herz}.

The dual of the Banach algebra $A_{p}(G)$ is the Banach space $PM_{p}(G)$ consisting of all limits of convolution operators associated to bounded measures \cite[Chapter 4]{Derigheti}.

The group $G$ is said to be amenable if, there exists an $m\in L^{\infty}(G)^{*}$ such that $m\geq 0$, $m(1)=1$ and $m(L_{x}f)=m(f)$ for each $x\in G$ and $f\in L^{\infty}(G)$ where $L_{x}f(y)=f(x^{-1}y)$ \cite[Definition 4.2]{Pier}.% (See also \cite{Derigheti} for more details of Figa-Talamanca Herz algebras and an amenable group $G$).

First we give the following lemma which is a generalization of Leptin-Herz Theorem (\cite[Theorem 10.4]{Pier}).
%======================================================================================================
\begin{lemma}\label{Lem: L-H}
 Let $G$ be a locally compact group and $1<p<\infty$. Then $A_{p}(G)$ has a b.$\Delta$-w.a.i if and only if $G$ is amenable.
\end{lemma}
\begin{proof}
Let $\{e_{\alpha}\}$ be a b.$\Delta$-w.a.i for $A_{p}(G)$ and $e\in A_{p}(G)^{**}$ be a $w^{*}$-cluster point of $\{e_{\alpha}\}$.

So, for each $\phi\in \Delta(A_{p}(G))=G$, we have
$$<e,\phi>=\lim_{\alpha}\phi(e_{\alpha})=1\cdot$$

Therefore, by \cite[Proposition 2.8]{Ulger} $G$ is weakly closed in $PM_{p}(G)=A_{p}(G)^{*}$. Now, by \cite[Corollary 2.8]{Chou} we conclude that $G$ is an amenable group.
\end{proof}
%====================================================================================
For a locally compact group $G$ let $L^{1}(G)$ be the group algebra of $G$ endowed with the norm $||.||_{1}$ and the convolution product as defined in \cite{Hewitt}.  By \cite[Theorem 23.7]{Hewitt2} we know that
\begin{equation*}
\Delta(L^{1}(G))=\{\phi_{\rho}; \rho\in \widehat{G}\},
\end{equation*}
where $\widehat{G}$ is the space of all continuous homomorphisms from $G$ into the circle group $\mathbb{T}$ and $\phi_{\rho}$ defined by
\begin{equation*}
\phi_{\rho}(h)=\int_{G}\overline{\rho(x)}h(x)dx\hspace{0.5cm}(h\in L^{1}(G))\cdot
\end{equation*}
%===========================================================================================================
\begin{theorem}\label{theorem} Let $G$ be a  locally compact  group.
 \begin{enumerate}
   \item If $G$ is an amenable group, then $L^{1}(G)$ is  $\Delta$-weak $\phi$-amenable for each $\phi\in \Delta(L^{1}(G))\cup\{0\}$.
   \item For $1<p<\infty$ and $\phi\in \Delta(A_{p}(G))\cup\{0\}$, $A_{p}(G)$ is $\Delta$-weak $\phi$-amenable   if and only if $G$ is an amenable group.
 \end{enumerate}
\end{theorem}
%%%%%%%%%%%%%%%%%%%%%%%%%%%%%%%%%%%%%%%%%%
%%%%%%%%%%%%%%%%%%%%%%%%%%%%%%%%%%%%%%%%%%
%Are we can remove the amenability of G?If we done this we construct a Banach algebra that has a B.D.A.I but does not have a B.A.I.
%%%%%%%%%%%%%%%%%%%%%%%%%%%%%%%%%%%%%%%%%%
%%%%%%%%%%%%%%%%%%%%%%%%%%%%%%%%%%%%%%%%%%
\begin{proof}
(1): It follows from  \cite[Corollary 2.4]{Monfared} that $L^{1}(G)$ is $\phi$-amenable for all $\phi\in \Delta(L^{1}(G))$. Also, we know that the group algebra has a bounded approximate identity and this completes the proof of (1) by using Theorem \ref{main}.

(2): If $G$ is an amenable group, then by Leptin-Herz's Theorem  we know that $A_{p}(G)$ has a bounded approximate identity. On the other hand, by  \cite[Corollary 2.4]{Monfared}, $A_{p}(G)$ is $\phi$-amenable for each $\phi \in \Delta(A_{p}(G))$. So, the result follows from Theorem \ref{main}.

Conversely, let $A_{p}(G)$ be $\Delta$-weak $\phi$-amenable. If $\phi=0$ the result follows from Lemma \ref{Lem: L-H}. If $\phi\in \Delta(A_{p}(G))$,  by Corollary \ref{corollary} we know that $A_{p}(G)$ has a b.$\Delta$.w.a.i. So, the result follows from Lemma \ref{Lem: L-H}.
\end{proof}

\section{Examples}
In this section we only give some illuminating examples.

The following example shows that the $\Delta$-weak 0-amenablity and  0-amenablity are different.
\begin{example} $A=l^{1}(\mathbb{Q}^{+})$ is $\Delta$-weak 0-amenable, but it is not 0-amenable. Because, $A$ has a b.$\Delta$-w.a.i, but it does not have any approximate identity \cite{Jones}.
\end{example}
The following three examples give Banach algebras which are $\Delta$-weak $\phi$-amenable, but they are not  $\phi$-amenable.
\begin{example} \label{example1}Let $X$ be a Banach space and take $\phi\in X^{*}\setminus \{0\}$ with $||\phi||\leq 1$. Define a product on $X$ by $ab=\phi(a)b$ for all $a, b\in X$. With this product $X$ is a Banach algebra which we denote it by $A_{\phi}(X)$. It is clear that $\Delta(A_{\phi}(X))=\{\phi\}$ and $A_{\phi}(X)$ is $\phi$-amenable if and only if dim$(X)=1$ \cite[Example 2.4]{Nasr}.

 Thus, if we take a Banach space $X$ with dim$(X)>1$ and a non-injective $\phi\in X^{*}$, then $A_{\phi}(X)$ is not $\phi$-amenable, but it is $\Delta$-weak $\phi$-amenable by  Proposition \ref{mainprop}.
 
On the other hand, let $x_{0}\in X$  be such that $\phi(x_{0})=1$. Then $x_{0}$ is a $\Delta$-weak identity for $A_{\phi}(X)$, but it is not an identity. Because, for  $0\neq a\in \textmd{ker}(\phi)$ we have $ax_{0}=\phi(a)x_{0}=0$. Therefore, $ax_{0}\neq a$. So, $x_{0}$ is not an identity. Moreover, it is clear that the $\Delta$-weak identity of the Banach algebra $A_{\phi}(X)$ is not unique, since each element $a_{0}\in X$ such that $\phi(a_{0})=1$, is a $\Delta$-weak identity for $A_{\phi}(X)$.

\end{example}
%==========================================================================
Let $A$ and $B$ be Banach algebras with $\Delta(B)\neq \emptyset$ and $\theta\in \Delta(B)$. The $\theta$-Lau product $A\times_{\theta}B$ is defined as the Cartesian product $A\times B$ with the following multiplication,
\begin{center}
  $(a,b)(a_{1},b_{1})=(aa_{1}+\theta(b)a_{1}+\theta(b_{1})a,bb_{1})\hspace{0.5cm}(a,a_{1}\in A, b,b_{1}\in B)\cdot$
\end{center}
With the $l^{1}$-norm and the above multiplication,  $A\times_{\theta}B$ is a Banach algebra. For see a full discussions on these kind of Banach algebras see \cite{Monfared3}.
\begin{example} Let $A=B=A_{\phi}(X)$ be the Banach algebra which defined in the previous example and dim$(X)>1$. Consider the $\phi$-Lau product $A_{\phi}(X)\times_{\phi}A_{\phi}(X)$. By   \cite[Proposition 2.8]{Monfared3} we know that $|\Delta(A_{\phi}(X)\times_{\phi}A_{\phi}(X))|=2$. Let $\Delta(A_{\phi}(X)\times_{\phi}A_{\phi}(X))=\{\Theta_{1},\Theta_{2}\}$. So, by proposition \ref{mainprop}, $A_{\phi}(X)\times_{\phi}A_{\phi}(X)$ is $\Delta$-weak $\Theta_{1}$-amenable and $\Delta$-weak $\Theta_{2}$-amenable. On the other hand, $A_{\phi}(X)$ is not $\phi$-amenable. Hence, by \cite[Lemma 6.8 (iii)]{Monfared2} there exists a character  $\Theta_{i}$, $i=1$ or 2  of $A_{\phi}(X)\times_{\phi}A_{\phi}(X)$ such that $A_{\phi}(X)\times_{\phi}A_{\phi}(X)$ is not $\Theta_{i}$-amenable.
\end{example}
%========================================================
\begin{example} Let $n\geq 2$ be an integer number and let $A$ be the Banach algebra of all upper-triangular $n\times n$ matrix over $\mathbb{C}$. We have $\Delta(A)=\{\phi_{1}, \phi_{2}, \phi_{3},..., \phi_{n}\}$ where
\begin{equation*}
\phi_{k}([a_{ij}])=a_{kk}\hspace{0.5cm}(k=1, 2, 3, ..., n)\cdot
\end{equation*}
Then for each $\phi_{k}$, $A$ is $\Delta$-weak $\phi_{k}$-amenable. To see this let  $e_{0}=[a_{ij}]$ be an element of $A$ such that
\begin{equation*}
a_{ij}=
\begin{cases}
1& i=j, i\neq k\\
0& i=j=k \\
0& i\neq j
\end{cases}\cdot
\end{equation*}
Obviously, $e_{0}$ is in $\ker(\phi_{k})$. But $\phi_{i}(e_{0})=1$ for all $i\neq k$. So, the result follows by using Theorem \ref{T: main charac}.

Also, $A$ is not $\phi_{k}$-amenable for each $k\geq 2$. Because, $\ker(\phi_{k})$ does not have a  right  identity. Therefore, by \cite[Proposition 2.2]{Kaniuth} $A$ is not $\phi_{k}$-amenable.
%this example obtained from example 6.5 of the paper, on character amenable Banach algebra by Sangani Monfared.
\end{example}
%===========================================================================

Some of Banach algebras satisfies both concepts of $\phi$-amenablity and $\Delta$-weak $\phi$-amenablity.
\begin{example}\label{ex1} If $A$ is a  $C^{*}$-algebra, then $A$ is  $\phi$-amenable and $\Delta$-weak $\phi$-amenable, because each $C^{*}$-algebra and their closed ideals have bounded approximate identity \cite[Theorem 3.2.21]{Dales}.
\end{example}
There exists a Banach algebra that is nighter  $\phi$-amenable nor $\Delta$-weak $\phi$-amenable as the next example shows. Also  the following example shows that the converse of Corollary \ref{Cor: 0-amen} is not valid in general.
\begin{example}\label{ex2} Let $A=C^{1}[0,1]$  be the Banach algebra consisting of all continuous function on $[0,1]$ with continuous derivation and norm $||f||_{1}=||f||_{\infty}+||f^{'}||_{\infty}$. We know that
\begin{center}
$\Delta(A)=\{\phi_{t}: \phi_{t}(f)=f(t) \textmd{ for each } t\in[0,1]\}\cdot$
\end{center}
By  \cite[Example 2.5(I)]{Kaniuth}, $A$ is not $\phi_{t}$-amenable for any $t\in[0,1]$.
Moreover,  there does not exist $t_{0}\in[0,1]$ such that $A$ is $\Delta$-weak $\phi_{t_{0}}$-amenable. To see this let $\{f_{n}\}$ be a b.$\Delta$-w.a.i for  $\ker(\phi_{t_{0}})$. So, it has the following properties
\begin{enumerate}
%  \item $\lim_{n}f_{n}(t_{0})=0$,
  \item $\lim_{n}f_{n}(t)=1\hspace{0.5cm}(t\in[0,1]\setminus\{t_{0}\})$,
  \item $\{||f^{'}_{n}||_{\infty}\} \textmd{ is bounded}\cdot$
\end{enumerate}
Hence, there exists a non-negative constant $M$ with $||f^{'}_{n}||_{\infty}=\sup_{t\in[0,1]}|f^{'}_{n}(t)|<M$ for all $n\in \mathbb{N}$. Hence, for positive integer $n_{0}$ we have
\begin{equation*}
\lim_{t\rightarrow t_{0}}|\frac{f_{n_{0}}(t)-f_{n_{0}}(t_{0})}{t-t_{0}}|=|f^{'}_{n_{0}}(t_{0})|<M\cdot
\end{equation*}
Therefore, there exists $\epsilon>0$ such that for each $t\in N(t_{0},\epsilon)=\{t: 0<|t-t_{0}|<\epsilon\}$ we have
\begin{equation*}
|f_{n_{0}}(t)-f_{n_{0}}(t_{0})|<M|t-t_{0}|\cdot
\end{equation*}
But the above relation is not valid in general, because the right hand side of the  inequality tends to zero as $t\rightarrow t_{0}$, but the left hand side do not.

Hence, $A$ is not $\Delta$-weak $\phi_{t}$-amenable for each $t\in[0,1]$.

Also, this Banach algebra is $\Delta$-weak 0-amenable, because the sequence $\{\frac{n-t^{n}}{n}\}$ is a bounded $\Delta$-weak approximate identity for $A$.
\end{example}
The converse of Theorem \ref{T: HP dense range} does not hold in general as the following example shows.
\begin{example} Let $A=C^{1}[0,1]$, $B=C[0,1]$ and $h:A\hookrightarrow B$ be the inclusion homomorphism. It is clear that $A$ is dense in $B$. By Examples \ref{ex1} and \ref{ex2} for each $t\in [0,1]$, $B$ is $\Delta$-weak $\phi_{t}$-amenable, but $A$ is not $\Delta$-weak $\phi_{t}$-amenable.
\end{example}
%%%%%%%%%%%%%%%%%%%%%%%%%%%%%%%%%%
%%%%%%%%%%%%%%%%%%%%%%%%%%%%%%%%%%
%further comments
%1-Please investigate the $\Delta$-weak $\phi$-amenability of $\teta$-Lau product.
%2-we can define the concept of n-delta-weak ap.id with use of n-homomorphisms. indeed replace each homomorphism with n-homomorphism.
%%%%%%%%%%%%%%%%%%%%%%%%%%%%%%%%%%%%%%
%%%%%%%%%%%%%%%%%%%%%%%%%%%%%%%%%%%
\bibliographystyle{amsplain}

\end{document}